\documentclass[10.75pt]{article}

\usepackage{amsthm}
\usepackage{times}
\usepackage[T1]{fontenc}
\usepackage{mathptmx}

\usepackage{amsmath, amssymb,color,float}
\usepackage{graphicx}
\usepackage{accents}
\usepackage{mathtools}
\usepackage{amsfonts}
\usepackage{times}
\usepackage{tabularx}
\usepackage[caption = false]{subfig} 
\usepackage{enumerate}
\usepackage{comment}
\usepackage{bm}
\usepackage[utf8]{inputenc}
\usepackage[english]{babel}
\usepackage{array}
\usepackage[a4paper, total={6in, 8in}]{geometry}
\usepackage{etoolbox}

\makeatletter
\patchcmd{\@makecaption}
  {\parbox}
  {\advance\@tempdima-\fontdimen2} % decrease the width!
  {}{}
\makeatother  

\usepackage{authblk}
\usepackage{natbib}
\newcommand{\Blp}{\Big(}
\newcommand{\Brp}{\Big)}

\newcommand{\rv}{\mathrm{Var}}
\newcommand{\E}{\mathrm{E}}

\newtheorem{theorem}{Theorem}

\newtheorem{corollary}{Corollary}

\theoremstyle{definition}
\newtheorem{example}{Example}
\newtheorem{remark}{Remark}

 \title{General Behaviour of P-Values Under the Null and Alternative}
\author[1,2]{Yanbo Tang  \thanks{Corresponding author, yanbo@utstat.toronto.edu, Supplementary Materials available on request.}}
\author[1]{ Radu Craiu }
\author[1,3]{Lei Sun}
\affil[1]{Department of Statistical Sciences, University of Toronto~~}
\affil[2]{Vector Institute}
\affil[3]{Division of Biostatistics, Dalla Lana School of Public Health, University of Toronto}

\begin{document}

\maketitle

\begin{abstract}
Hypothesis testing results often rely on simple, yet important assumptions about the behaviour of the distribution of $p$-values under the null and \textit{alternative}. 
We examine tests for one dimensional parameters of interest that converge to a normal distribution, possibly in the presence of many nuisance parameters, and characterize the distribution of the $p$-values using techniques from the higher order asymptotics literature. 
We show that commonly held beliefs regarding the distribution of $p$-values are misleading when the variance or location of the test statistic are not well-calibrated or when the higher order cumulants of the test statistic are not negligible.   
Corrected tests are proposed and are shown to perform better than their first order counterparts in certain settings. 
\end{abstract}

%\keywords{$p$-values; \and multiple hypothesis testing; \and  FDR; \and FWER; \and Higher Order Asymptotics }

\section{Introduction}
Statistical evidence against a hypothesis often relies on the asymptotic normality of a test statistic, as in the case of the commonly used Wald or score tests. 
Many authors ignore the asymptotic nature of the argument and assume that in finite samples the distribution of the test statistic is indeed normal.
This perfunctory approach generates misleading beliefs about the $p$-value distribution, such as i) the distribution of the $p$-values under the null is exactly uniform or that ii) the cumulative distribution function (henceforth, cdf) of the $p$-value under the alternative is concave. 
However, there are important exceptions from these rules, e.g. discrete tests are not normally distributed in any finite sample settings, so that the distribution of the $p$-values under the null is certainly not uniform. 
Similarly, it is not obvious that the cdf is concave under the alternative as we will illustrate with some examples.
Testing procedures aimed at controlling the family-wise error rate (FWER) or the false discovery rate  (FDR, see \cite{BH}) typically assume that i) or ii) holds.
In \cite{Cao2013}, the authors examine the optimality of FDR control procedures when i) or ii) are violated and provide alternative conditions to maintain said optimality.
Clearly, a more precise characterization of the $p$-value distribution that accounts for the approximation error is pivotal in controlling the occurrence of false discoveries.

 Complicating matters even further, the issue of calibrating the location and variance of the test statistic is often overlooked, particularly under the alternative.
 Under the alternative the test statistic can be improperly re-scaled since often the variance of the test statistic is obtained under the null.
 While under the null, the test statistic may not have zero mean and may also not be correctly standardized, thus making the standard Gaussian approximation suspect.
 The problem of biases in the variance and expectation is aggravated in the presence of a large number of nuisance parameters. 
For instance, while it has been demonstrated in \cite{DiCicio} that the profile score statistic has a location and variance bias under the null,  in Section \ref{subsec:score} we show that the variance bias can persist under the alternative.

%Assuming that the test statistic is exactly normally, \cite{hung} studied the $p$-value distribution, characterizing the density of the $p$-value distribution under the alternative.
%\textcolor{blue}{However, there is no work characterizing the behaviour of $p$-values accounting for the approximations involved on the test statistic, especially under the alternative, which \cite{Cao2013} showed is important to fully characterize optimality of the false discovery control (FDR) procedure.}

%The approximate nature of most $p$-values used is practice has not often been acknowledged in the literature, and is somewhat surprising, as $p$-values are a critical component of statistical hypothesis testing. 
%\textcolor{blue}{Furthermore assumptions on the distribution of the $p$-values under the alternative are generally more implicit than those made for the uniformity of the null for most testing procedures. \cite{Cao2013} explored this in detail and obtained that for optimally for false discovery rate control it is needed for the ratio of of the density of the $p$-values to be monotone (add more). Also FDR control procedures such as \cite{storey2002} requires that the distribution under the CDF alternative is concave  }

 These concerns motivate us to perform a systematic study of the $p$-value distribution in the presence of information or location biases under the null and alternative, while accounting for the approximation error resulting from the use of asymptotic arguments.
We explore how certain asymptotically non-vanishing and vanishing biases in the variance and location of the test statistic can occur in finite samples, violating the assumptions generally placed on the null and alternative distributions of the $p$-values. 
We study both continuous and discrete distributions supported on lattices.
In doing so we include all approximation errors, including those induced by discreteness, to fully characterize the behaviour of the distributions of $p$-values under the null and alternative.
This work extends the results of \cite{hung}, who studied the distribution of $p$-values under alternative assuming the test statistic is normally distributed, to a broader framework. 
% Our examination reveals that unexpected behaviour of $p$-values  can emerge using these finite sample approximations, and that sample size is a relative notion not wholly encapsulated by the number of observations, $n$. We then proceed to show how these unexpected behaviour can affect common procedures for controlling the family-wise error rate, (FWER).
%This can occur with some types of score tests where the variance is assumed known and in high dimensional settings when the number of parameters are comparable with the number of data-points such as in \cite{Sur2019}, where it was shown that the variance of the likelihood ratios statistic is underestimated in high dimensional logistic regression. 
 
We focus on univariate test statistics for a one dimensional parameter of interest based on sums of independent random variables, possibly in the presence of a large number of nuisance parameters. 
These types of test statistics are commonly used to infer the significance of individual coefficients in most regression models.
The results of the paper are in the same vein as those found in \cite{hall2013bootstrap} and \cite[\S~3]{kolassa1994series}, whose objective was the coverage properties of confidence intervals. 
We expand their results to the $p$-values, motivated by the multitude of scientific investigations that rely on the $p$-value distribution rather than confidence intervals.

\begin{figure}[ht]
\centering
\subfloat{\includegraphics[width = 2in]{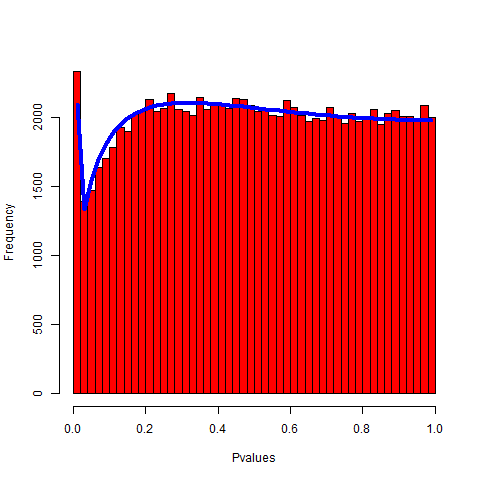}} 
\subfloat{\includegraphics[width = 2in]{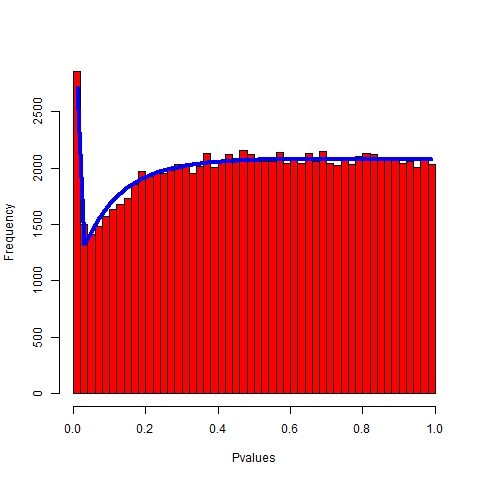}}\\
\subfloat{\includegraphics[width = 2in]{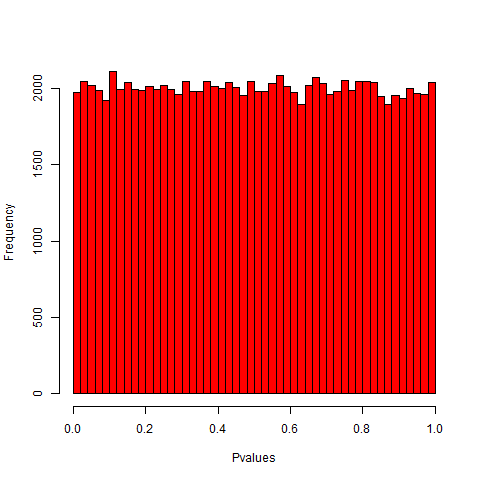}}
\subfloat{\includegraphics[width = 2in]{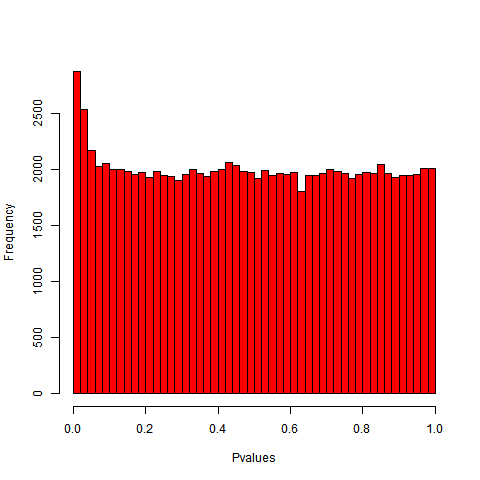}}
\caption{Distribution of $p$-values under $H_0$ and $H_1$ in Example 1. {\it Upper left}: $p$-values obtained under $H_0$ from the normal approximation. {\it Upper right}: $p$-values obtained under $H_1$ from the normal approximation. {\it Lower left and right}:  corrected $p$-values for the null and alternative, respectively, using the saddlepoint approximation that is introduced in Section 3. The number of samples is 750. The upper left panel clearly does not exhibit uniformity and upper right panel's distribution does not appear to have a concave cdf. We plot the theoretical prediction from Theorem 1 in blue for the upper left and upper right panels.}
\label{fig:test}
\end{figure}

We begin with a simple example illustrating how the standard assumptions on the null and alternative distributions of the $p$-values can be violated in practice.
%that sample size alone cannot provide exhaustive information about the  convergence rate of a statistical test, and motivate the need for a more precise characterization of the $p$-value distribution. 
\begin{example}
We wish to test the null hypothesis $H_0: \beta = 0.01 $ against the alternative $H_1: \beta = 0.01/1.05$, where $\beta$ is the rate parameter of a gamma distribution, based on 750 observations $x_1, \cdots, x_{750}$, assuming that the shape parameter is known to be $\alpha = 0.01$.  From the central limit theorem, we know that the test statistic 
\begin{align*}
S_n = \sqrt{n} \left(\frac{\bar{X} - 1}{\sigma} \right) \rightarrow N(0, 1),
\end{align*}
so we are able to obtain a two-sided $p$-value based on the standard normal distribution, where $\sigma = Var(X_1)$. 
We plot the histograms of the $p$-values obtained under the null and alternative in Figure \ref{fig:test}. 
The plots are obtained by simulation using 100,000 replications.  

We see on Figure \ref{fig:test} that the distribution of the $p$-values obtained from the simulations does not adhere to its expected behaviour under the null or the alternative. 
The upper left plot in Figure \ref{fig:test} shows a marked departure from the $U(0,1)$ distribution expected from the null. Thus, a typical rejection rule which assumes uniformity of the $p$-value distribution under the null will not provide type I error control for certain choices of $\alpha$.
For example, if we desire a $10^{-4}$ significance level, we obtain a type I error approximately equal to $1.579*10^{-3}$, which is fifteen times higher than the nominal level.
Under a local alternative, the upper right plot in Figure 1 shows that the $p$-value distribution may not be stochastically smaller than a $U(0,1)$. 
The resulting lack of concavity of the distribution $p$-value under the alternative can violate the typical assumption that the false negative rate is strictly decrease and the FDR is increasing in the nominal control level $\alpha$ in the multiple testing setting; see \cite{Cao2013}.  
Note that the cause for this poor calibration is not the low sample size.
\end{example}

\begin{figure}[ht]
\centering
\includegraphics[width=8cm, height = 5cm]{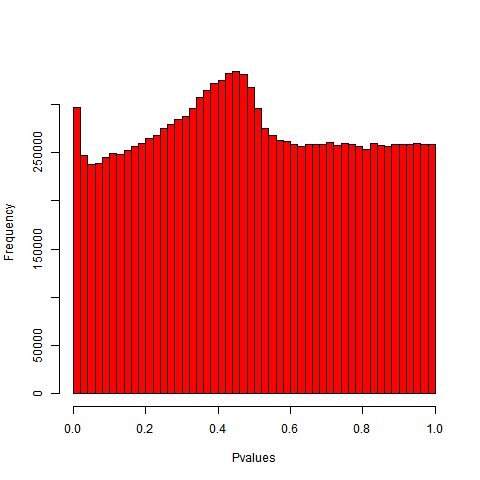}
\caption{Empirical p-value distribution based on a genome wide association study ($n=45,637$) of lung cancer. 
}
\label{fig:lung}
\end{figure}

To assure the reader that the above example is not a singular aberration, we present in Figure \ref{fig:lung} the histogram of over 13 million $p$-values from the genome-wide association study of lung cancer generated from the UK Biobank data \citep{biobank}. These $p$-values are produced by the Neale Lab \citep{neale}, based off  45,637 participants and 13,791,467 SNPs; SNPs with minor allele frequency less than 0.1\% and INFO scores less than 0.8 were excluded from the analysis.
We note that the histogram exhibits a similar behaviour to the one seen in Example 1, i.e.,  the distribution of $p$-values exhibits a secondary mode that is far from zero. 

\begin{figure}[H]
\centering
\subfloat[Shape 1]{\includegraphics[width = 2in]{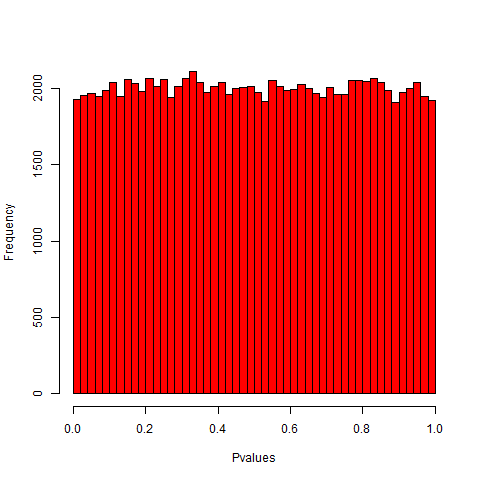}} 
\subfloat[Shape 2]{\includegraphics[width = 2in]{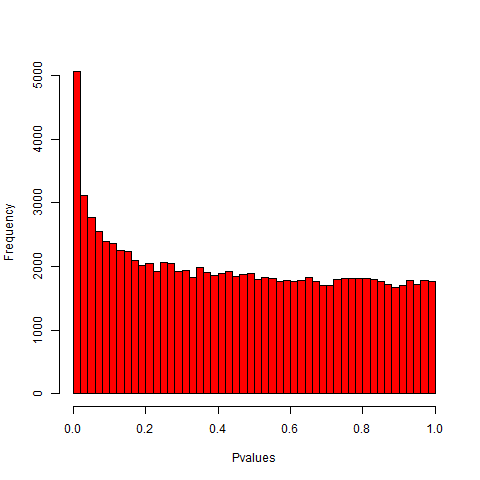}}\\
\subfloat[Shape 3]{\includegraphics[width = 2in]{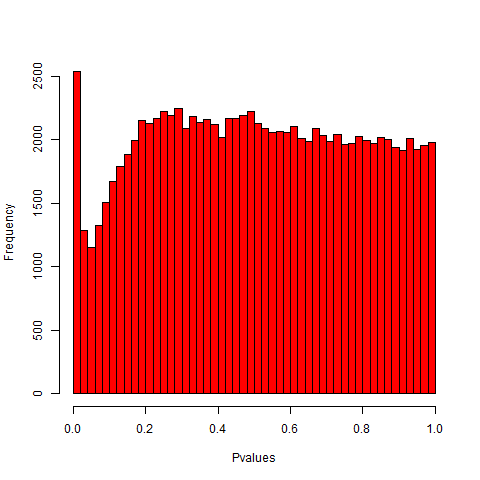}}
\subfloat[Shape 4]{\includegraphics[width = 2in]{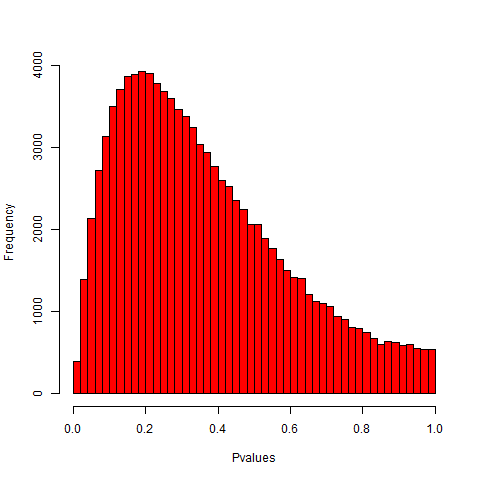}} \caption{General chart of the behaviour of $p$-values under the null and alternative for a two tailed test. Shapes 1 to 3 were obtained from simulating from the null and alternative from Example 1 with different parameters $\alpha$ and $\beta$. Shape 4 was obtained when using a misspecified variance, as detailed in Section 2.2. }\label{shapes}
\end{figure}

% Within this paper we will try to first study the general shapes that the $p$-value density can take under the null and alternative, taking into account the fact that in certain scenarios we are quite far from the limiting normal distribution. 
Figure \ref{shapes}  briefly summarizes the shapes that the density of the $p$-value distribution might take for two-tailed tests, based on the results in Theorems 1 and 2 that are introduced in Section 2. The descriptions in Table 1 %\ref{tb:shapes} 
verbalise the various mathematical conditions that can lead to the four shapes in Figure \ref{shapes}. 
In practice it is possible to have combinations of the shapes listed in Figure \ref{shapes}, as the observed test statistics may not be identically distributed and can be drawn from a mixture of the null and alternatives hypotheses.

\begin{table}[h]
\caption{Description of the test statistic's characteristics and the resulting shapes (as shown in Figure \ref{shapes}) of the $p$-value distribution under the null or alternative hypothesis.}
\fbox{% 
\begin{tabular}{ | m{3em} | m{6cm}| m{6cm} |} 
 \hline
 Shape & Null & Alternative  \\ 
 \hline
 1 & The typical uniform shape. & Possible if effect size is small.  \\ 
 \hline
 2 & Possible if variance is misspecified, underestimated. & Typical behaviour.   \\
 \hline
 3 & Possible if test statistic has large higher order cumulants, see Example 1. & Possible if effect size is small and the higher order cumulants are large, see Example 1. \\
 \hline
 4 & Possible if variance is misspecified, overestimated, see Example 3. & Possible if variance is misspecified, overestimated and the effect size is small, see Example 3.  \\
 \hline
 \end{tabular}}
\end{table}

Section 2 contains the main theoretical results of this paper, Theorems 1 and 2, which characterize the distribution of $p$-values under the null and alternative. 
Section 2.1 examines the $p$-value distribution resulting from the score test, while Section 2.2 studies specific examples.
Section 3 provides numerical results and considers some remedies aimed at calibrating the $p$-value distribution. 
Section 4 closes the paper with a discussion of the implications of our results and some recommendations to practitioners.
 
\section{Distribution of $p$-values under Non-Normality}
%This section contains the main theoretical result of the paper which describes the general behaviour of $p$-values. 
All theoretical details and proofs, as well as a brief introduction of the concepts needed for the proof of Theorems 1 and 2, are deferred to the Supplementary Materials. 
We consider the case where the test statistic, $S_n$, can be discrete and may also have a non-zero mean under the null and a non-unit variance under the null or alternative.
We assume that $\psi$ is a one dimensional parameter of interest and $\lambda$ is a vector of nuisance parameters.
Without loss of generality, let the statistic $S_n$ either be used to test the null hypothesis $H_0: \psi = \psi_0$ for a two-sided test or $H_0: \psi \geq \psi_0 $ for a one-sided test.
All results are given in terms of the cdf. 
%as the density with respect to the Lebesgue measure for discrete variables is not a well-defined notion. 
Theorems 1 and 2 deal with the case where $S_n$'s distribution is continuous and discrete respectively. 
% After the presentation of the Main theorems we will provide numerical examples in \S2.1 and some practical considerations for multiple hypothesis testing corrections in \S2.2.

We first consider the case where the statistic $S_n$ admits a density. 
We typically assume that $S_n$ has been appropriately calibrated such that $\E[S_n]= 0$ under the null, and $\E[S_n] \neq 0$ under the alternative hypothesis.
We let $p(S_n)$ denote the $p$-value obtained from the test statistic $S_n$.
However, as discussed in the introduction, the mean of $S_n$ may not be exactly $0$ under the null due to a location bias.
We also would expect that the variance of the test statistic should be 1 under the null and alternative, which may not be the case for all test statistics however; see Example 3.
The location bias complicates the precise determination of whether $S_n$'s distribution should be considered under  the null or the alternative.
However, note that Theorem 1 statement is applicable under both the null and alternative, since its conclusion depends only on the expectation, variance, and the other cumulants of the statistic $S_n$, regardless of the true hypothesis.
We first introduce some notations:
\begin{itemize}
\item[(i)] $Z_p$ is the $p$-th quantile of a standard normal distribution. 
\item[(ii)] $\rho_{n,i}$ is the $i$-th order standardized cumulant of $S_n$, and $\rho_n$ is a vector containing all cumulants.
\item[(iii)] $\phi(x)$ is the standard normal density.
\end{itemize}

\begin{theorem}\label{th:cont_approx}
Let $X_1,\ldots, X_n$ be a sequence of continuous, independent random variables. Set $S_n=  \sqrt{n} (\bar{X}_n - a_n)/b_n$ where $\bar{X}_n =  n^{-1}\sum_{i = 1}^n X_i$, and let $\{a_n\}_{n\ge 0}$, $\{b_n\}_{n\ge 0}$ be two sequences of real numbers. Let $\E[S_n]= \mu_n$, ${\rv}(S_n) = v_n^2$, and $\rho_n$ denote the cumulants of $(S_n - \mu_n)/v_n$.  Then the CDF of the one-sided $p$-value is
\begin{align}
    \mathbb{P} (p(S_n) < t)  = \Phi\left(\frac{Z_t - \mu_n}{v_n}\right) -E_2\left(\frac{Z_t - \mu_n}{v_n}, \rho_{n} \right) + O\left(n^{-3/2}\right),
\end{align}
and the CDF of the two-sided $p$-value is:
\begin{align}
    \mathbb{P} (p(S_n) < t)  &= 1 + \Phi\left(\frac{Z_{t/2} - \mu_n}{v_n} \right) - \Phi\left(\frac{-Z_{t/2} - \mu_n}{v_n}\right)\nonumber \\
    &+E_2\left(\frac{Z_{t/2} - \mu_n}{v_n} , \rho_n\right) - E_2\left( \frac{-Z_{t/2} - \mu_n}{v_n}, \rho_n\right) + O\left(n^{-3/2} \right), \label{eq:cont_one}
\end{align}
where, 
\begin{align}
E_2(t, \rho_n) =   -\phi(t)\Big\lbrace  \frac{\rho_{n,3} H_2(t)}{6} + \frac{\rho_{n,4} H_3(t)}{24} + \frac{(\rho_{n,3} )^2 H_5(t)}{72} \Big\rbrace, \label{eq:cont_two}
\end{align}
and $H_j(t)$ denotes the j-th Hermite polynomial.
\end{theorem}

\begin{remark}
The $j$-th order Hermite polynomial is a polynomial of $j$-th degree defined through the differentiation of a standard normal density. A table of the Hermite polynomials is given in the Supplementary Materials. 
\end{remark}

\begin{remark}
Should an approximation to the probability density of the $p$-value distribution be desired, it can be obtained from differentiating Equations (\ref{eq:cont_one}) and (\ref{eq:cont_two}).
\end{remark}

In general $E_2 = O(1/n^{1/2})$, however in the
the case that $\mu_n = 0$,  $E_2(Z_{t/2}/ v_n , \rho_n) - E_2( -Z_{t/2}/v_n, \rho_n) = O(1/n)$ for two-sided tests due to cancellations which occur in the difference of the odd Hermite polynomials.
We refer to terms in $E_2(t, \rho_n)$ as the higher order terms. 
Therefore, supposing $\mu_n = 0$ under the null, meaning the sequence $a_n = \mathbb{E}[\bar{X}_n]$, we obtain the following corollary: 

\begin{corollary}
Assume the setting and notation from Theorem \ref{th:cont_approx} and suppose that under the null we have $\E[S_n]= 0$, and $\rv(S_n) = 1$. The CDF of the distribution of the $p$-values for a one-sided test under the null is
\begin{align}
    \mathbb{P} \left(p(S_n) < t \right)  = t + O\left( n^{-1/2}\right),
\end{align}
and the CDF of the distribution of the $p$-values for two-sided test under the null is
\begin{align}
    \mathbb{P} \left(p(S_n) < t \right)  &=  t + O\left(n^{-1} \right). 
\end{align}
\end{corollary}

Corollary 1 shows that the two-sided test is preferable unless there is a scientific motivation for using the one-sided test. 
The case when $S_n$ has a discrete distribution supported on a lattice is covered in Theorem \ref{th-discrete}. 
\begin{theorem}
\label{th-discrete}
Let $X_1, \cdots, X_n $  be a sequence of independent discrete random variables where $X_i$ has mean $m_i$. Suppose that $X_i - m_i$ is supported on a lattice of the form $c + j\cdot d$, for $j \in \mathbb{Z}$ and for all $1\le i \le n$. Assume $d$ is the largest number for which this property holds. 
Set $S_n=  \sqrt{n} (\bar{X}_n - a_n)/b_n$,  where $\bar{X}_n =  n^{-1}\sum_{i = 1}^n X_i$, $\E[S_n]= \mu_n$,  $\rv(S_n) = v_n^2$,  $\rho_n$ as the cumulants of $(S_n - \mu_n)/v_n$ and $d_n = d   /(\sqrt{n}v_n b_n) $.

Then the CDF of the one-sided $p$-value is
\begin{align*}
    \mathbb{P} (p(S_n) < t) &=  \Phi\left( \frac{Z_t - \mu_n}{v_n} \right) +E_2 \left( \frac{ Z_t - \mu_n}{v_n}, \rho_n \right)  
     + C_2\left(\frac{ Z_t - \mu_n}{v_n}, \rho_n \right) + O\left(n^{-3/2}\right),
\end{align*}
and the CDF of the two-sided $p$-value is
\begin{align*}
    \mathbb{P} (p(S_n) < t)  &=1 + \Phi\left( \frac{Z_{t/2} - \mu_n}{v_n} \right) - \Phi\left(\frac{-Z_{t/2} - \mu_n}{v_n} \right) 
    +E_2\left(\frac{Z_{t/2} - \mu_n}{v_n}, \rho_n \right) \\
    &- E_2 \left(\frac{-Z_{t/2} - \mu_n}{v_n}, \rho_n \right)   +C_2 \left(\frac{Z_{t/2} - \mu_n}{v_n}, \rho_n \right) 
    - C_2\left(\frac{-Z_{t/2} - \mu_n}{v_n}, \rho_n \right)  +O \left(n^{-3/2} \right),
\end{align*}
where,
\begin{align*}
C_2(t, \rho_n) =  - d_n Q_1\Blp \frac{ t - \sqrt{n}c}{d_n}\Brp \Blp 1 + \frac{\rho_{n,3} H_3(t)}{6} \Brp + \frac{ d_n^2}{2} Q_2\Blp \frac{ t - \sqrt{n}c}{d_n}\Brp,
\end{align*}
and $Q_j(t)$ are periodic polynomials with a period of $1$. On $[0,1)$, they are defined by
\begin{align*}
Q_1(t) = t - \frac{1}{2}, \quad Q_2(t) = t^2 - 2t +\frac{1}{6},
\end{align*}
and $E_2(t, \rho_n)$ and $H_j(t)$ are defined as in Theorem 1.
\end{theorem}

% \begin{remark}
% In the case where we are considering a two tailed test and $\mu_n$ is 0, the odd polynomials cancels out.  Resulting in a simpler expression for the $p$-value distribution under the null. Working through the algebra one can determine that if $\rho_3^n \neq 0$ or $\rho_4^n >0$ then all of the observed approximated $p$-values close to 0 will be too small compared to the true $p$-value.
% \end{remark}

\begin{corollary}
Assume the setting and notation from Theorem \ref{th-discrete} and suppose that under the null $\E[S_n]= 0$, and $Var(S_n) = 1$. Then the $p$-values obtained from one or two-sided tests satisfy
\begin{align}
    \mathbb{P} (p(S_n) < t)  = t + O\left(n^{-1/2}\right).
\end{align}
\end{corollary}
Note that the convergence  is slower by a factor of $n^{-1/2}$ compared to the continuous case for a two-sided test.  This is due to the jumps in the CDF which are of order $O(n^{-1/2})$.

\begin{remark}
Under the alternative, the $p$-value distribution depends on the effect size, $\mu_n$, as well as the magnitude of the higher order cumulants.
However, for large values of $\mu_n$ the impact of the higher order terms will be negligible, as $E_2$ is a product of an exponential function and a polynomial function which decays to 0 asymptotically in $\mu_n$. We explore this further in Example 4.
\end{remark}

\begin{remark}
When performing multiple hypothesis testing corrections, the $p$-values of interest are often extremely small. 
Therefore from Corollary 1 and 2, we see that a large amount of samples is needed to guarantee the level of accuracy required since the approximation error is additive.
\end{remark}

%\begin{remark}
%It is possible to approximate a general statistic $S_n = \sqrt{n}(\hat{\psi} - \theta_0)/\hat\sigma$ that converges to a normal distribution, the approximation take the same form as proposition 1,  although the polynomials involved in the approximation may change as mentioned in section 2.1. The convergence to normality for Wald type tests can also be characterized by the sizes of the higher order cumulants. This is discussed in the Supplementary Materials.
%\end{remark}

\subsection{An Application of the Main Theorems: The Score Test} \label{subsec:score}
 We examine the broadly used score test statistic, also known as the Rao statistic. 
 The popularity of the score statistic is due to its computational efficiency and ease of implementation.  
In the presence of nuisance parameters, the score statistic is defined through the profile likelihood. 
Suppose that the observations $y_i$'s are independent, and denote the log-likelihood function by $l(\psi, \lambda; Y) = \sum_{i = 1}^n l(\psi, \lambda; y_i) $ then
\[ l_\text{pro}(\psi) = \sup_{\lambda} l(\psi, \lambda; Y) = l(\psi, \hat\lambda_\psi; Y), \] 
where $\hat\lambda_\psi$ denotes the constrained maximum likelihood estimator. The score statistic is defined as:
\[ S_n(\psi_0) = \frac{l_{\text{pro}}^\prime(\psi_0)}{ \lbrace -l^{\prime\prime}_{\text{pro} } (\psi_0) \rbrace^{1/2}}= \sum_{i = 1}^n \frac{ \frac{d}{d\psi} l(\psi, \hat\lambda_\psi; y_i)}{ \left\lbrace-\sum_{i = 1}^n \frac{d^2}{d\psi^2} l(\psi, \hat\lambda_\psi; y_i)\right\rbrace^{1/2}} \xrightarrow{D} N(0,1),  \]
under the usual regularity assumptions. 
Due to the form of $S_n$, we may apply Theorem 1 or 2.

The presence of nuisance parameters induces a bias in the mean and variance of the score statistic; see \cite{profile_bias} and \cite{DiCicio}. 
Thus, it is not the case that the mean of the score statistic is 0 and the variance is 1 under the null, as the profile likelihood does not behave like a genuine likelihood and does not satisfy the Bartlett identities. 
In general this problem is compounded if the number of nuisance parameters is increased, as we illustrate below. 

We only discuss the location bias, since the formulas for the information or variance bias are much more involved and compromise the simplicity of the arguments. 
From \cite{profile_bias}, the bias of the profile score under the null is:
\begin{align}
\mathbb{E}\lbrace l^\prime_{\text{pro}}(\psi_0) \rbrace &= \alpha_n + O\left( n^{-1} \right),
\end{align}
where the term $\alpha_n = O(1)$. The form of $\alpha_n$ is given in the Supplementary Materials.

% \begin{align}
% \mathbb{E}\lbrace l^\prime_{\text{pro}}(\psi_0) \rbrace &= -\frac{1}{2}\sum_{i,j = 1}^{p}(\kappa_{\psi, ij} - \sum_{k,l = 1}^{p}\kappa_{\psi,k} \kappa^{k,l} \kappa_{l. ij})k^{i,k} - \sum_{i,j = 1}^{p} \frac{1}{2}(\kappa_{\psi, i,j} -  \sum_{k,l = 1}^{p} \kappa_{\psi, k}\kappa^{k,l}\kappa_{l, i,j} )\kappa^{i,j} \nonumber\\
% &+ O\left( n^{-1} \right),
% \end{align}

% where in the above $\kappa$ with lowered indices denotes the expectation of likelihood derivatives, a comma indicates a multiplication of derivatives, $i,j,k,l$ are reserved for the nuisance parameter $\lambda$, for example:
% \begin{align*}
%     \kappa_{ij, k} = \mathbb{E}\left[ \frac{\partial^2 l(\psi, \lambda; Y)}{\partial\lambda_i\partial\lambda_j} \frac{\partial l(\psi, \lambda; Y)}{\partial\lambda_k}  \right] .
% \end{align*} 
% We use $\kappa$ with raised indices to denote the inverse of matrices, for example:
% \begin{align*}
%      \kappa^{i,j} = \left( \mathbb{E}\left[ \left\lbrace \frac{\partial l(\psi, \lambda; Y)}{\partial\lambda} \frac{\partial l(\psi, \lambda; Y)}{\partial\lambda^\top} \right\rbrace^{-1} \right] \right)_{ij}.
% \end{align*}

To estimate the effect of the dimension of the nuisance parameter on the size of the bias, we use a similar argument as \cite{laplace}, in which they count the number of nested summations that depends on $k$, the number of parameters in the model, to estimate the rate of growth of a function in $k$. 
From the expression of $\alpha_n$ given in the Supplementary Materials, we obtain at most 4 nested summations which depends on $k$ therefore the bias of the profile score is of order $O(k^4)$ in the worst case scenario. 
The rather large location bias can be impactful as it may induce a perceived significance when $k$ is large, an example using Weibull regression is given in Section 3.3. 
A similar argument can be applied to the information bias; see \cite{DiCicio} for a comprehensive discussion on the form of these biases. 

The information bias for the score statistic can also be highly influential under the alternative.
In that case, the expected value of the score statistic is non-zero, which is desirable, but the variance of the statistic $S_n$ can be either over- or underestimated. 
Since the true parameter value is not $\psi_0$, there is no guarantee that ${l_{\text{pro}}^{\prime\prime}(\psi_0)}$ gives the correct standardization. 
If the estimated variance is larger than the true variance of the score, then it is possible to obtain Shape 4 in Figure \ref{shapes}, which violates the concavity assumption for $p$-value distribution's CDF under the alternative. 
Further, if we assume that under the null the $p$-value distribution is uniform, then this also violates the monotonicity assumption required by \cite{Cao2013} for the optimality of FDR control.
Example \ref{ex:score_glm} illustrates this phenomenon using the score test in a generalised linear model.

\begin{example} \label{ex:score_glm}
Assume the following regression model based on the linear exponential family where the density of the observations $y_1, \cdots, y_n$ are independent and follows
\[ h(y_i| \beta, X_i) = \exp\lbrace a(X_i\beta) y_i+ b(X_i\beta) + D(y_i)  \rbrace, \]
where $X_i$ is a vector of covariates associated with each $y_i$, and $\beta = (\beta_0, \beta_1, \cdots, \beta_k)$ is a vector of regression coefficients. 
Let $f(X\beta) = E[y|X]$ denote the mean function. 
\cite{score_reg} studied the score statistic for testing the global null $\beta_1 = \beta_2= \cdots = \beta_k = 0$ and linked the resulting statistic to linear regression. 
A similar analysis can be performed for different hypothesis, such as inference for a parameter of interest in the presence of nuisance parameters to produce a more general result whose derivation is consigned to the Supplementary Materials. The resulting score statistic takes the form 
\[S_n = \lbrace y - f(X\hat\beta_{\text{null}}) \rbrace^\top W X \left\lbrace X^\top D X \right\rbrace^{-1} X^\top W \lbrace y - f(X\hat\beta_{\text{null}} ) \rbrace \xrightarrow[]{D} \chi^2_{q}, \]
where $f(X\hat\beta_{\text{null}})$ is a vector whose $i$-th entry is $f(X_i\hat\beta_{\text{null}})$, $q$ is the number of constraints in the null hypothesis, and $\hat\beta_{\text{null}}$ denotes the constrained maximum likelihood estimate under the null. W and D are square diagonal matrices of dimension $n$ whose entries are $[W]_{ii} = a^\prime(X_i\hat\beta_{\text{null}})$ and $[D]_{ii} = a^\prime(X_i\hat\beta_{\text{null}}) f^\prime(X_i\hat\beta_{\text{null}}) $ for $i = 1, \dots, n$. Using a suitable change of variable, the statistic $S_n$ can be related to weighted linear regression.

In the common case where we wish to test for $\beta_j = 0$, the score statistic can be re-written in the form:
\[S_n = [ (X^\top D X)^{-1} ]_{jj}^{1/2} \sum_{i = 1}^n a^\prime(X_i \hat\beta_{\text{null}}) x_{ij} \lbrace y_i - f(X_i \hat\beta_{\text{null}}) \rbrace  \xrightarrow{D} N(0,1),  \]
under the null. Under the alternative we may write 
\[S_n = (1 + c_n )\tilde{S}_n + d_n,  \]
where $\tilde{S}_n$ converges in distribution to a standard normal. The scaling factor $c_n = O(1)$ is an information bias and $d_n$ plays the role of the effect size and will increase to infinity as the number of samples increases. 
However, if $\beta_j \approx 0 $ then $d_n$ can be quite small, meaning that the effect of the scaling factor $c_n$ can be consequential. 
For an example of this see Example 3, where the effect size is not large enough to offset the scaling factor. 

\end{example}

\begin{remark}
Although the likelihood ratio statistic can be written as a summation of independent random variables, the limiting distribution of the likelihood ratio test is a gamma random variable, therefore Theorem 1 or 2 are not directly applicable.
It may be possible to modify the baseline density used in the Edgeworth expansions to obtain a result based on Laguerre polynomials.
This can also be useful when examining the asymptotic behaviour of test statistics for testing vector parameters of interest, as these test statistics often have a gamma distributed limiting distribution. 
\end{remark}

\begin{remark}
 The bias issue discussed within this section is also present for the Wald test statistics, even if it can not be represented by a summation of independent random variables. 
It is rarely the case that the maximum likelihood estimate is unbiased, and the same applies for the estimate of the variance of the maximum likelihood estimate. 
Generally the problem worsens as the number of nuisance parameters increases.
\end{remark}

\subsection{Numerical Examples of Application of the Main Theorems}
We illustrate the results of the main results with some numerical examples to demonstrate how various problems in the distributions of the $p$-value can occur. 
We first examine a discrete case where the statistic $S_n$ does not admit a density. 
We note that when the $E_2$ term is negligible, our results on the distribution of the $p$-values coincide with those obtained by \cite{hung} when the exact normality of the test statistic holds.  
 
On the contrary, when the additional terms are not negligible or the variance is incorrectly specified,
the behaviour of the distribution of $p$-values can be quite different.  
The exact size of the difference depends on the behaviour of the Hermite polynomials, the higher order cumulants and the variance. 
We consider the following examples in order to illustrate some of the ramifications.

\begin{example}\label{ex:linkage}
Consider a simple linkage analysis of sibling pairs who share the same trait of interest, a common problem in statistical genetics.  
The underlying principle is that genes that
are responsible for the trait are expected to be over-shared between relatives, while the null hypothesis states that the trait similarity does not impact allele sharing, i.e. independence between the trait and gene.
The problematic distribution of $p$-values in this example is caused by the discrete nature of the problem along with a misspecified variance under the alternative.

Since the offsprings are from the same parents, under the null we would expect the number of shared alleles to be either 0, 1 or 2 with probability $\theta_{null} = (p_0, p_1, p_2) = (0.25,0.5,0.25)$ based on Mendel's first law of segregation.
However, under the alternative we can expect the sharing level to be higher than expected. 

Assume that we have $n$ affected sibling pairs.
Let $x_i$ be the number alleles shared amongst the $i$-th affected sibling pair. 
Then under the null $E[x_i] = 1$, $Var[x_i] = 0.5$ and we let $y_i = (x_i - 1)/\sqrt{0.5}$.
We consider the following well known non-parametric linkage test 
\begin{align*}
    S_n = \sum_{i = 1}^n \frac{y_i}{\sqrt{n}}  \sim N(0,1),
\end{align*}
see \cite{laird2010fundamentals}. The above can be compared to a score test as only the information under the null was used.

Under the alternative, the distribution of the test can be misspecified, since a different distribution of allele sharing  will yield a different variance. 
Consider the simple example when the distribution of the numbers of shared alleles follows a multinomial distribution with $\theta_{alt1} = (0.09,0.8,0.11)$.
The variance of this distribution is $0.4 < 0.5$. Yet another alternative in which $\theta_{alt2} = (0.29,0.4,0.31)$ yields the variance  $0.6 > 0.5$; in both case there is oversharing.
We include visualizations of  the $p$-value distribution under the two alternatives in Figure \ref{fig:linkage}. Theorem 2 is used to produce the approximation given by the blue curve, due to the discrete nature of the test statistic.

\begin{figure}
    \centering
     \subfloat[Score test, $\theta = \theta_{alt1}$]{\includegraphics[width = 2in]{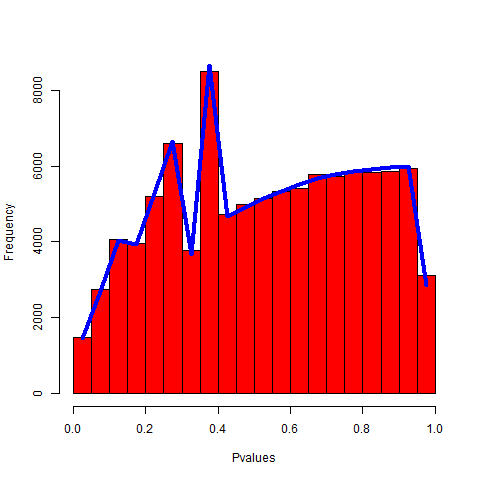}}
    \subfloat[Score test, $\theta = \theta_{alt2}$ ]{\includegraphics[width = 2in]{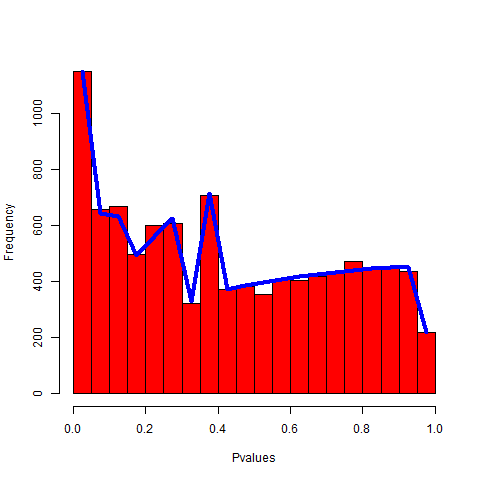}}
     \\
    \subfloat[Wald test,  $\theta = \theta_{alt1}$]{\includegraphics[width = 2in]{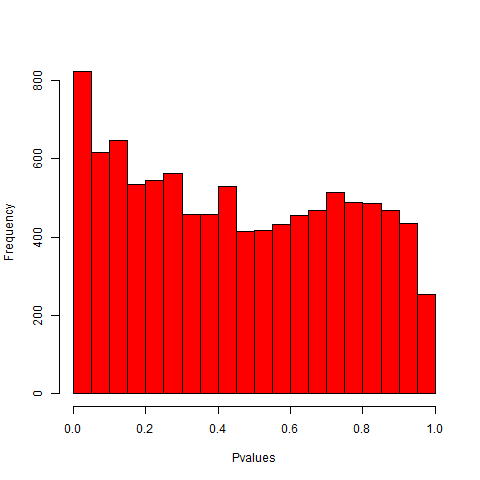}} 
    \subfloat[Wald test, $\theta = \theta_{alt2}$ ]{\includegraphics[width = 2in]{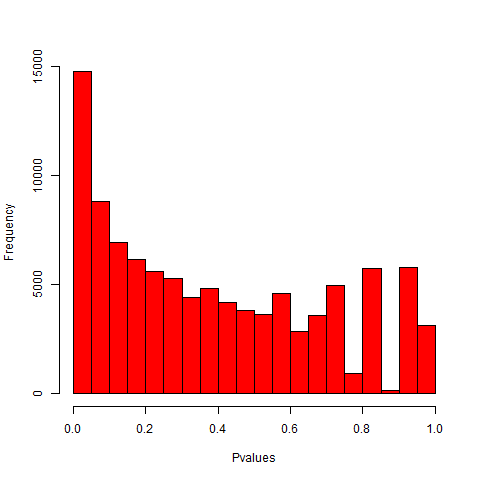}}
    \caption{Plots for Example 3, examining the behaviour of the $p$-value distribution for non-parametric linkage analysis for the score test (upper panel) and the Wald test (lower panel). The simulation is performed with $n =400$, and $100,000$ replications. Samples of sibling pairs are generated from a multinomial distribution with $\theta_{alt1} = (0.09,0.8, 0.11) $ for the two plots on the left panel, and $\theta_{alt2} = (0.29,0.4,0.31)$ for the two plots on the right panel. For the score test, both histograms have spikes due to the discrete nature of the problem. The discrete version of the Edgeworth approximation, plotted in blue, is used as the test statistic is supported on a lattice. The histograms of the $p$-values obtained form the Wald test look much better than their score test counterparts. }
    \label{fig:linkage}
\end{figure}
\end{example}

In this case the problem can be resolved by considering a Wald type test where the variance is calculated from the maximum likelihood estimate $\hat{\theta} = ( \#(x_i = 0)/n, \#(x_i = 1)/n, \#(x_i = 2)/n  )$, and use:
\[S_n^\prime = \sum_{i = 1}^n \frac{x_i - 1}{\sqrt{n\widehat{\text{var}}(x_i)}}.
\] 
We plot the results of applying the Wald test in Figure \ref{fig:linkage}.
The solution is quite simple in this case, but in more complex models it is more computationally expensive to calculate the variance estimate under the alternative. \\

\noindent\textbf{Example 1 revisited.} The abnormal distribution of $p$-values in this scenario is caused by a large numerical value of $\rho_{n,3}$ and $\rho_{n,4}$.
Going back to Example 1, we look at the theoretically predicted behaviour of the $p$-values under the null and alternative.
Figure \ref{fig:test} shows the histograms of the empirical $p$-values obtained by simulation versus the theoretical prediction given in Theorem 1, shown as the blue curve. 
Without accounting for the higher order terms in the expansion we would have expected the null distribution to be uniform, however, using Theorem 1, we obtain a much more accurate description of the $p$-value distribution. 
%Under the null, we have successfully predicted the presence of a large concentration of $p$-values near 0, which should signal that we cannot naively assume that the test statistic is normal as this would cause a highly inflated type I error for small values of $\alpha$.  
In the bottom panel of Figure \ref{fig:test} we also show a corrected version of the $p$-values approximation using the the saddlepoint approximation which will be introduced in Section 3.

The estimation of small p-values based on the standard normal approximation can be drastically optimistic. We report in Table \ref{tb:deltas} 
the differences between the exact and the approximate $p$-value obtained from Example 1 for the 5 smallest $p$-values.
The smallest $p$-values from the normal approximation are not on the same scale as the exact $p$-values, the smallest approximate $p$-value being five-fold times smaller than its exact counterpart.
In contrast, the $p$-values produced by the saddlepoint approximation are very close to the exact ones.
\begin{table}[b]
\caption{ \label{tb:deltas} Table of $p$-values obtained from Example 1 under the null. The exact $p$-values are obtained from the density of the gamma distribution, the approximate $p$-values are obtained from the normal approximation.}
\centering
\begin{tabular}{rrrrr}
  \hline
 ID & rank & $p$-value exact & $p$-value approx. & $p$-val saddlepoint \\ 
  \hline
60326 &   1 & 1.04E-05 & 1.04E-10 & 1.04E-05 \\ 
  91132 &   2 & 1.46E-05 & 3.06E-10 & 1.47E-05 \\ 
  83407 &   3 & 2.12E-05 & 9.66E-10 & 2.12E-05 \\ 
  97470 &   4 & 3.31E-05 & 3.75E-09 & 3.32E-05 \\ 
  2573 &   5 & 3.80E-05 & 5.66E-09 & 3.81E-05 \\ 
  \hline
\end{tabular}
\end{table}

%The approximation error in the tail of the $p$-value distribution is especially relevant to multiple hypothesis testing, as often the $p$-values of interest are on the scale of $10^{-6}$ or possibly even smaller.

\begin{example}
 We examine the influence of the effect size $\mu_n$ on the distribution of the $p$-values under the alternative using the same set-up as in Example 1. In our simulations we increase the effect size $\mu_n$ by changing the value of $\beta$, while keeping $\alpha$ fixed. The results are displayed in Figure \ref{fig:altnernaives}.

 \begin{figure}[H]
\centering
\subfloat[Subfigure 1 list of figures text][Small effect size  ]{\includegraphics[width = 1.75in]{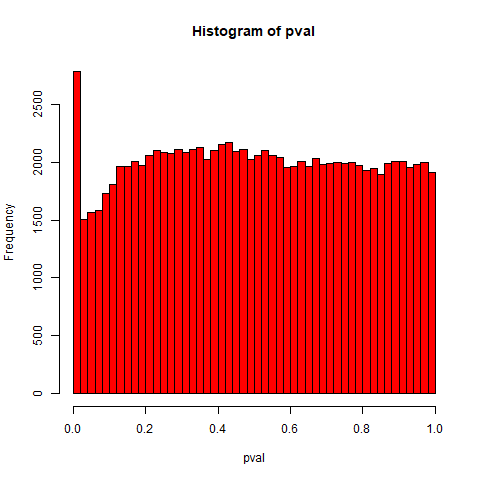}} 
\subfloat[Subfigure 2 list of figures text][Medium effect size ]{\includegraphics[width = 1.75in]{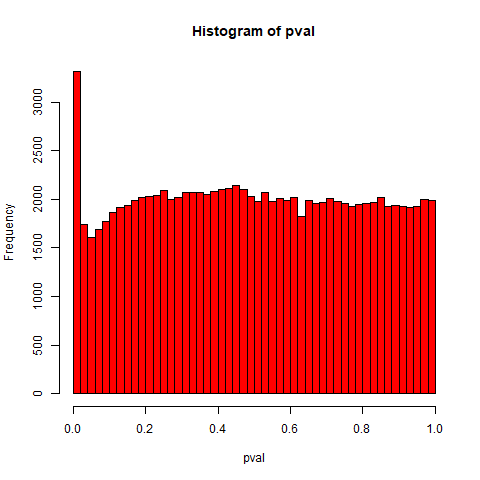}}
\subfloat[Subfigure 1 list of figures text][Large effect size ]{\includegraphics[width = 1.75in]{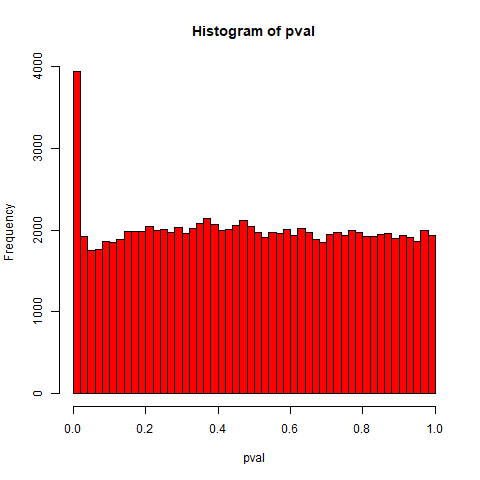}} \\
\subfloat[Subfigure 1 list of figures text][Small effect size ]{\includegraphics[width = 1.75in]{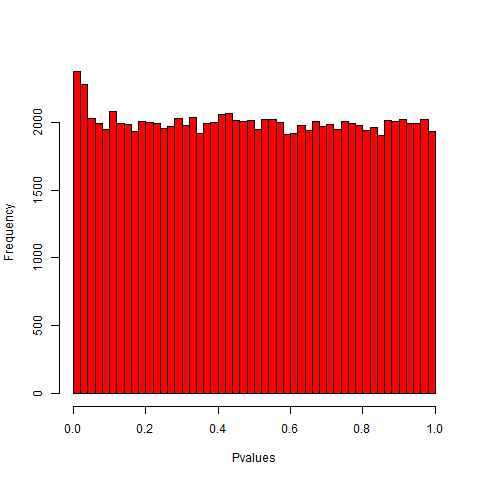}} 
\subfloat[Subfigure 2 list of figures text][Medium effect size ]{\includegraphics[width = 1.75in]{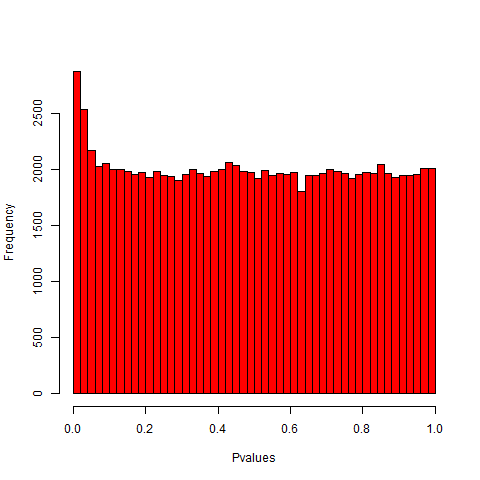}}
\subfloat[Subfigure 1 list of figures text][Large effect size ]{\includegraphics[width = 1.75in]{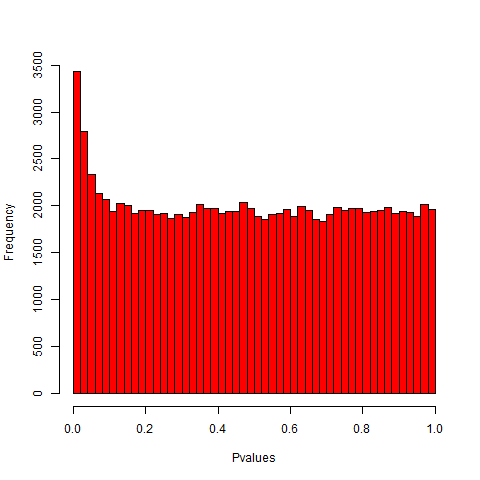}}
\caption{Distribution of the approximated $p$-values (top panel) and the corrected $p$-values (bottom panel), under three different alternatives with $\alpha_1 =\alpha_2 = \alpha_3 =0.01$ and $\beta_2 = 0.01/1.025$, $\beta_1 = 0.01/1.05$ and $\beta_3 = 0.01/1.1$, from left to right.
 By increasing the effect size, the approximate $p$-values starts to behave in an expected manner. 
 While the corrected $p$-values obtained by using the saddlepoint approximation is well behaved for all effect sizes.}
\label{fig:altnernaives}
\end{figure}

As discussed in Remark 3, for large effect sizes $\mu_n$, the distribution of $p$-values generated from the test statistic follows the expected trend, where there is a concentration of $p$-values around $0$ and the density decreases in a monotone fashion to 1. 
Conversely, should $\mu_n$ be small then the behaviour under the alternative can be quite different from what we would expect, as illustrated by the top-right plot in Figure \ref{fig:altnernaives}.
\end{example}

\section{Additional Examples and Possible Remedies}
We provide additional examples of problematic $p$-value distributions, and we explore some possible remedies based on high order asymptotics. 
We also provide additional examples of problematic $p$-value distributions.

A commonly used tool for higher order asymptotics is the saddlepoint approximation, which is a density approximation that can be integrated to obtain tail probabilities, e.g. $p$-values. 
For a good survey of the saddlepoint approximation and its applications in statistics, we refer the reader to \cite{reid1988} or for a more technical reference, we suggest \cite{jensen1995saddlepoint} or \cite{kolassa1994series}.

The saddlepoint approximation can be most easily obtained for a sum or average of independent random variables, $X_1, \dots, X_n$. The density approximation then results in an approximation of the cumulative distribution through a tail integration argument,
\begin{align}
    P(\bar{X} < s) = \Phi(r_s)\lbrace 1 + O(n^{-1}) \rbrace, \label{accuracy_saddle}
\end{align}
where $r_s$ is a quantity constructed from the saddlepoint and the cumulants of the distribution of the $X_i$'s. This can be used for conditional inference in generalized linear models by approximating the distribution of the sufficient statistics in a exponential family model; see \cite{davison}.

% \begin{remark}
% An improved approximation can be obtain by re-normalizing the saddlepoint approximation such that it integrates to 1. 
% The re-normalization improves the accuracy of the approximation by a further factor of $1/\sqrt{n}$.
% \end{remark}
Another more broadly applicable tail approximation is the normal approximation to the $r^\star$ statistic \citep{cox1994inference}, which is obtained by adding a correction factor to $r$, the likelihood root. It can be used in regression settings for inference on a scalar parameter of interest.
Let $r =\text{sign}(\hat\psi - \psi_0) [2 \lbrace l_{\text{pro}}(\hat\psi) - l_{\text{pro}}(\psi_0) \rbrace]^{1/2}$ denote the likelihood root, and in what follows the quantity $Q$ varies depending on the model.
\begin{align*}
    P(r < r_{obs}) = \Phi\left\lbrace r_{obs} + \frac{1}{r_{obs}} \log\left(\frac{Q}{r_{obs}}\right)\right\rbrace \left\lbrace 1 +O\left( n^{-3/2} \right) \right\rbrace,
\end{align*}
where $r_{obs}$ is the value of the likelihood root statistic based on the observed data.
Using the above, we also obtain an improved approximation to the true distribution of the likelihood root. 
For a discussion of $r^\star$ see \cite{reid_wald}.

The proposed methods require two model fits, one under the alternative and one under the null in order to obtain $r$, contrary to the score test. 
The methods listed here are by no means comprehensive since there are a variety of other candidates which may be of use, such as the often applied Firth correction \citep{firth} or other forms of bias correction obtainable by adjusting the score equation \citep{kosmidis2020mean}.

\subsection{The Gamma example}
We apply the saddlepoint approximation to Example 1 and display the results in Figure \ref{fig:test}. 
Considering the null $H_0: \alpha = \beta = 0.01$ (the two plots on the left panel), there is a spike around 0 for $p$-values obtained using the CLT (top left plot).
In contrast, we see a marked improvement of the overall
behaviour of the $p$-value distribution after the proposed correction (bottom left plot). 

\subsection{Logistic Regression in Genetic Association Studies}
We apply the normal approximation to $r^\star$ to a simulated genome-wide association study to further illustrate the practical use of the proposed correction. 
We consider a logistic regression model that links the probability of an individual suffering from a disease to that individual's single nucleotide polymorphism (SNP), a genetic ordinal variable coded as 0, 1 or 2, and other covariates such as age and sex. 
Formally, let the disease status of the individual be $Y_i$, which is either $0$ (individual is healthy) or $1$ (individual is sick) and $\pi_i = E[Y_i]$ denote the probability of individual $i$ having the disease and let $X_{i, s}$ denote the genetic covariate of interest of the $i$-th individual, while $X_{i,j}, j = 1, 2$, $j\ne s$ are the other covariates. 
 The regression model is: 
\begin{align*}
\text{logit}(\pi_i) = X^i_s \beta_s + \sum_{j = 1}^2 X^i_j \beta_{j} + \beta_0.
\end{align*}
We consider the difficult case where the disease is uncommon in the population and the SNPs of interest are rare, i.e. most observed values of $X_{i,s}$ are 0. 
It is known that in this situation the single-SNP test performs poorly, and pooled analyses of multiple SNPs have been proposed \citep{pooled}. 
However for the purpose of this study, we assume that the individual SNPs are of interest.

\begin{figure}[ht]
\centering
\subfloat{\includegraphics[width=2.75in]{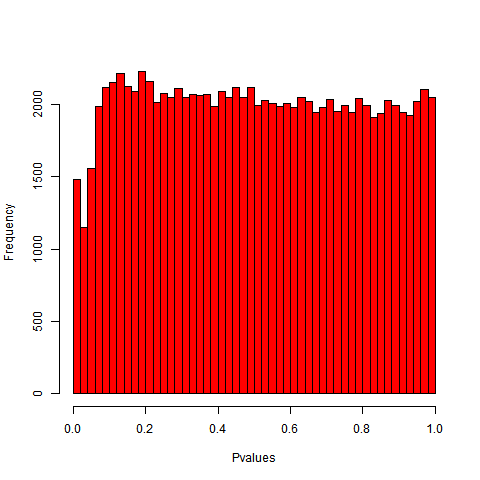}}
\subfloat{\includegraphics[width=2.75in]{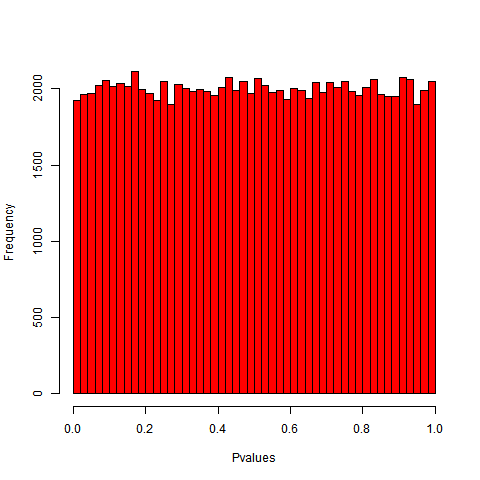}}
\caption{Empirical distribution of the null $p$-values from a logistic regression association study of SNPs with with low minor allele frequency and a low number of diseased individuals. 
Left histogram displays the $p$-values from the Wald test under the null. The right histogram displays the $p$-values histogram obtained from $r^\star$ under the null.  }
\label{fig:6}
\end{figure}

We consider a simulated example to demonstrate the effectiveness of the correction. 
We generate a sample of 3,000 individuals, their genetic variable $X_s$ are simulated from a $Binomial(2, 0.025)$, a binary variable $X_1$ from a $Binomial(1,0.5)$ and finally $X_2$ from a $N(20, 1)$. 
We let $\beta_0 = -3.5$, $\beta_s = 0$, $\beta_1 = 0.02$ and $\beta_2 = 0.02$. 
With this set of parameters we would expect on average $ \approx 4.6\%$ of the cohort to be in the diseased group, based on the expected value of the covariates, i.e. approximately 137 participants with $Y_i = 1$. 
For each replication of the simulation, we re-generate the labels from the logistic model. 
Figure \ref{fig:6} shows that the correction works well under the null. 

This example suggests that the usefulness of the proposed higher order corrections is not limited to small sample scenarios, as note by \cite{zhou2018efficiently} who used the saddlepoint approximation in case control studies with extreme sample imbalance. 
Naively we would expect that with 3,000 participants, of which 137 are in the diseased group, the Wald test should behave correctly. However, the skewed distribution of the SNP values severely reduces the accuracy of the test. 
The use of $r^\star$ corrects the distribution of the $p$-values as shown in Figure \ref{fig:6} (right plot) where the distribution of the $p$-values under the null ($\beta_s=0$) is approximately Unif(0, 1) as expected. 

In the example above it is clear that even though we have 3,000 individuals, of which 137 are affected by the disease, the standard approximation performs very poorly. 
This seems to suggest that in our particular example, the effective sample size is lower than 137 for the diseased group. 
Next we consider a simple regression with a single genetic covariate in order to illustrate the loss in information resulting from the sparsity of the minor allele. 
We use the available Fisher information about the parameter of interest as a measure of effective sample size. 

The standard deviation of the parameter of interest obtained from the inverse information matrix is
\begin{align*}
 var(\hat\beta_s) &= \frac{\sum_{i = 1}^{n} \hat{P_i} (1 - \hat{P_i})}{\sum_{i = 1}^{n} \hat{P_i} (1 - \hat{P_i})\sum_{i = 1}^{n} x_i^2 \hat{P_i} (1 - \hat{P_i}) - (\sum_{i = 1}^{n} x_i \hat{P_i} (1 - \hat{P_i}) )^2},\\
 &\approx \frac{1}{\sum_{i = 1}^{n} x_{i,s} \hat{P_i} (1 - \hat{P_i})},
\end{align*}
where $\hat{P}_i$ is the predicted probability of an individual being diseased and the approximation is valid under the assumption that the allele frequency is low enough such that we observe very few 1's and almost no 2's. 
The information about the parameter $\beta_s$ is increasing in terms of $x_i \hat{P_i}(1 - \hat{P_i})$. 
It is apparent that the rate of increase in information is limited by the sparsity of the rare allele. 
In order to have more information about the parameter, we would need to observe more individuals who have the rare allele, i.e.  $X_i \ne 0$.

\subsection{Logistic Regression - Data from the 1000 Genome Project}
We consider an additional logistic regression example as this type of model is broadly used in statistical genetics. 
Using phase 3 data from the 1000 Genome project \citep{1000genome}, we construct an artificial observational study in order to study how these approximations behave on real genome-wide genetic data. %examples of GWAS studies. 
In our simulations, we take the 2504 individuals within the database and assign the $i$-th individual a label of $0$ or $1$ based on the following logistic model, where $\pi_i = P(Y_i = 1)$:
\[ \text{logit}(\pi_i) =  \sum_{j = 1}^4 X^i_j \beta_{j}  + \beta_\text{Sex}*I( \text{Sex}_i = \text{male}) + \beta_0, \]
where $\text{Sex}_i$ is the biological sex of the $i$-th individual. Four other covariates are included, where $X^i_j$ are independent for all $i, j$ and follow a standard normal distribution.
The model coefficients are set to 
\[(\beta_0, \beta_1, \beta_2, \beta_3, \beta_4, \beta_{\text{Sex}}) = (-3.25, 0.025, -0.025, 0.025, -0.03, 0.1).\]
Once we assign a label to the $i$-th individual we keep it fixed throughout the simulation.
We then fit a logistic model using the SNPs for which the minor allele frequency is at least $1\%$ on chromosome $10$, and ethnicity as additional covariates.
We use the Wald test, and $r^\star$, but do not consider the cases where perfect separation occurs, as both methods considered here cannot deal with this issue. 

We plot some of the results for the Wald test and $r^\star$. 
We focus on rare variants with MAF $\leq 2.5\%$ and semi-common variants with $2.5\% <$ MAF $\leq 10\%$, as the remaining common variants are expected to behave well.  
In total $160,580$ SNPs fall into the rare variant category while $176,350$ SNPs fall into the semi-common variant category.

 \begin{figure}[h]
\centering
\subfloat[Subfigure 1 list of figures text][Rare variants, Wald]{\includegraphics[width =2in]{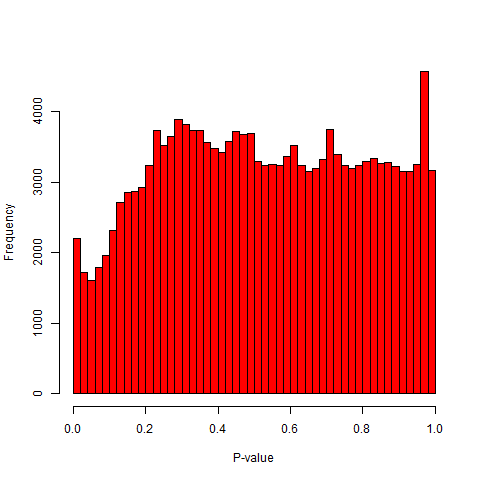}} 
\subfloat[Subfigure 1 list of figures text][Rare variants, $r^\star$]{\includegraphics[width = 2in]{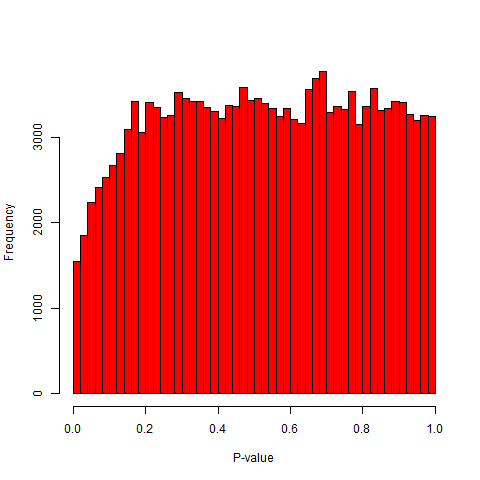}}\\
\subfloat[Subfigure 1 list of figures text][Semi-common variants, Wald ]{\includegraphics[width = 2in]{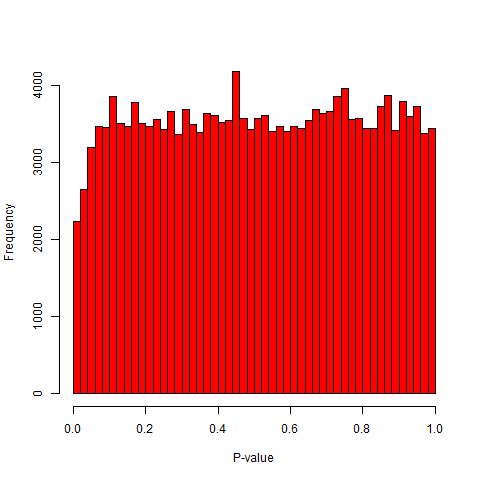}}
\subfloat[Subfigure 1 list of figures text][Semi-common variants, $r^\star$ ]{\includegraphics[width = 2in]{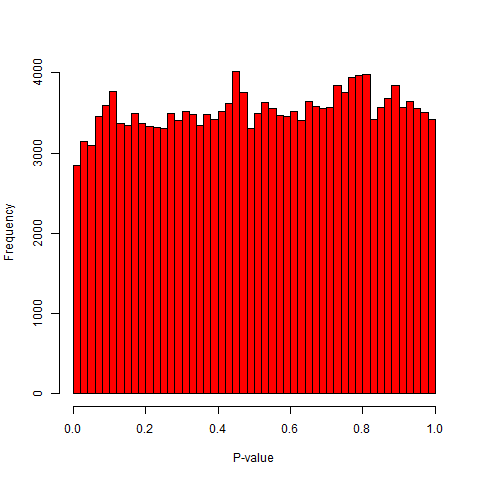}}
\caption{Distribution of $p$-values for Wald test and $r^\star$. The null distribution was simulated by using sex and four other randomly generated covariates. We fit a logistic regression model using SNPs from chromosome 10, with $160,580$ being rare ($\text{MAF}\leq 2.5\%$ ) and $176,350$ semi-common variants ($2.5\% < \text{MAF}\leq  10\%$ ).}
\label{altnernaives}
\end{figure}

As expected, the two tests behave better for semi-common SNPs than rare SNPs (bottom vs. top panel of Figure 7), producing $p$-values that more closely follow the Unif(0,1) distribution. Among the two tests, the proposed $r^\star$ method clearly out-performs the traditional Wald test. 
However, this application also points out the limitation of $r^\star$ as the correction for rare variants is not sufficient (top right plot), and further improvement of the method in this case is of future interest.

\subsection{Weibull survival regression}
Consider an example where there is a large number of nuisance parameters, leading to an inconsistent estimate of the variance.
We examine a Weibull survival regression model in which all of the regression coefficients, except the intercept, are set to $0$ by simulating $y_i \sim \text{Weibull}(1,2)$, independently of any covariate. 
We set the number of observations, $n$ to 200 and the number of covariates to $50$, and generated the covariates as IID standard Gaussian, and test for whether the first (non-intercept) regression coefficient is 0. 
We perform 10,000 replications and plot the histogram of the $p$-values, and compare the Wald test to the $r^\star$ correction. 

\begin{figure}[H]
\centering
\subfloat{\includegraphics[width=0.35\linewidth]{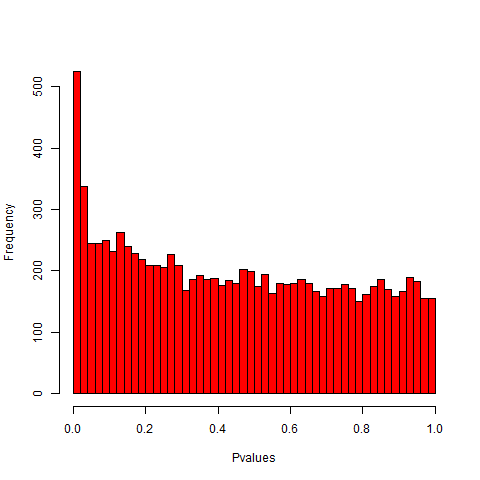}}
\subfloat{  \includegraphics[width=0.35\linewidth]{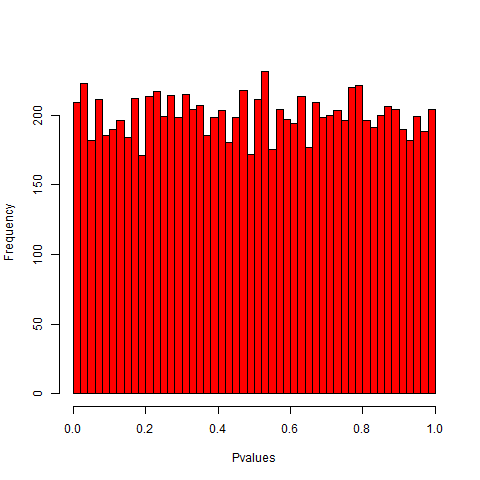}}
\caption{On the left, a histogram of the $p$-values produced by the Wald test for $\beta_1 = 0$ under the null with $n = 200$ and $p = 50$ with no censoring. On the right, a histogram of the $p$-value obtained from the $r^\star$ correction. 10,000 replications were performed.}
\label{regression}
\end{figure}

In Figure \ref{regression} we see a high concentration of $p$-values around $0$ for the Wald test, leading to increased I error.
The corrective procedure brings the distribution under the null much closer to uniformity.
We see that naively adding more and more information into the model while trying to perform inference on a one dimensional parameter of interest is problematic as it creates a perceived significance of the parameter of interest under the null.   

% \begin{remark}
% We briefly covered two classes of commonly used models in statistical applications, we refer the reader to \cite{2007applied} for a discussion of other classes of models where higher order asymptotic can be applied. 
% There is a package "likelihoodAsy", \cite{rstar} implemented in R that performs inference for an arbitrary likelihood function. 
% \end{remark}

\section{Discussion and Conclusion}
We characterize the distribution of $p$-values when the test statistic is not well approximated by a normal distribution by using additional information contained in the higher order cumulants of the distribution of the test statistic.
We also demonstrate that there are issues beyond failure to converge to normality in the that the expectation and variance of the test statistics can be misspecified, and these issues can persist even in large sample settings.
In doing so we have extended the previous work done by \cite{hung} to greater generality, examining the score test in exponential models in the presence of nuisance parameters.
We also examine some possible remedies for making the $p$-value distribution adhere more closely to their usual required behaviour such as uniformity under the null or concavity of the CDF under the alternative.
These assumptions are very important to justify the usage of current FWER and FDR procedures.
The proposed remedies may not solve all problems
relating to the $p$-value distribution in the finite sample settings, but they do at least partially correct some of the flaws.

We suggest the use of the proposed saddlepoint approximation or the normal approximation to $r^\star$ in practice, because a) the exact distribution of a test statistic is often unknown, b) the usual CLT approximation may not be adequate, and c) the high order methods are easy to implement. 
This will ensure a closer adherence to the assumptions usually needed to conduct corrective procedures used in FWER control or FDR control.

\section*{Acknowledgement}
The first author would like to thank Nancy Reid, Michele Lambardi di San Miniato and Arvind Shrivats for the help and support they provided. We also thank the Natural Sciences and Engineering Research Council, the Vector Institute and the Ontario government for their funding and support.

\newpage

\bibliography{biblio}

\end{document}